\documentclass[12pt]{amsart}
\usepackage{graphicx}
\usepackage{amssymb}
\usepackage{amsmath}
\usepackage{amsthm}
\usepackage{epstopdf}
\newcommand{\dif}{\mathrm{d}}

\begin{document}

\title{The Geodesics of Rolling Ball Systems}

\author{Daniel R. Cole}
\address{SUNY Maritime College}
\email{dcole@sunymaritime.edu}


\date{August 27, 2012}

\begin{abstract} A well-known and interesting family of sub-Riemannian space are the systems involving two balls rolling against each other without slipping or twisting.  In this note, we show how the sub-Riemannian geodesics of these space, when the two balls are embedded in $\mathbb{R}^3 \times \mathbb{R}^3$, are horizontal curves on the intersections of these balls with Euclidean $5$-planes.
\end{abstract}

\maketitle

\section{Introduction}

Rolling ball systems, in which two spheres with different radii roll against each other without slipping or twisting, are particularly interesting examples of sub-Riemannian manifolds.  E. Cartan in \cite{Cartan5manifolds} was the first to study these configuration spaces, albeit in a different form.  After the re-emergence of sub-Riemannian geometry as a topic of interest in the late twentieth century, rolling ball systems and sub-Riemannian manifolds with growth vector $(2,3,5)$ were studied by R. Bryant and L. Hsu in \cite{BryantRigidity}.  R. Montgomery's invaluable monograph on sub-Riemannian geodesics, \cite{geodesictour}, included a section on rolling ball systems, and an open question posed by him in that section about the relationship between the symmetries of rolling ball systems and the exceptional Lie group $G_2$ sparked a great deal of research over the past decade.  See \cite{ZelenkoRolling}, \cite{AgrachevRolling}, \cite{MontgomeryRolling}, and most recently \cite{BaezRolling} for more details.

Surprisingly, with all of this interest in rolling ball systems, we have yet to see an explicit description of the sub-Riemannian geodesics of a rolling ball system beyond the singular geodesics first discovered by A. Agrachev and A. Sarychev in \cite{AgrachevSingular}.  It turns out that, like the Heisenberg group, the sub-Riemannian geodesics of a rolling ball system have a very simple characterization.  Specifically, if we consider each sphere as embedded in $\mathbb{R}^3$ and centered at the origin, and the point of contact between the two spheres as an ordered pair in $\mathbb{R}^3 \times \mathbb{R}^3$, then for any given sub-Riemmanian geodesic of our rolling ball system, there exists a $5$-plane in $\mathbb{R}^3 \times \mathbb{R}^3$ such that all points of contact on the geodesic lie on this $5$-plane.  Since the rolling ball system comes equipped with a two-dimensional horizontal distribution, for almost all $5$-planes passing through a point, a unique horizontal curve is determined by the $5$-plane.  Thus this characterization of the sub-Riemannian geodesics of our system is quite specific.  

The proof below is quite simple, and should be understandable for an advanced undergraduate student.  Essentially, we choose five arbitrary points of contact along our candidate geodesic $\gamma$ and bend the curve at these points.  We then show that the coordinates of those points of contact in $\mathbb{R}^3 \times \mathbb{R}^3$ must be linearly dependent in order for $\gamma$ to be a length minimizing geodesic.  We note that a similar proof, using only three points, can be used to characterize the sub-Riemmanian geodesics of the Heisenberg group (and thus offers yet another proof of the isoperimetric problem in the plane).

\section{Results}

To begin, let $M_1$ and $M_2$ be the two spheres in our rolling ball system.  Without loss of generality, we may assume that $M_1$ has radius $1$.  Let $r$ be the radius of $M_2$.  Each sphere sits in its own copy of Euclidean space centered at the origin.  Therefore every point in $M_1$ can be represented by a unit vector $\mathbf{u}$, and every point in $M_2$ can be represented by $r\bf{v}$, a unit vector multiplied by $r$.  A point of contact between $M_1$ and $M_2$ consists of an ordered pair $\left( \mathbf{u}, r \mathbf{v} \right)$ of points in $M_1$ and $M_2$, respectively.  In general, a state in our rolling ball system can be represented by a point of contact and two unit vectors in $\mathbb{R}^3$, $\mathbf{a}$ and $\mathbf{b}$, where $\mathbf{u} \cdot \mathbf{a} = 0$ and $\mathbf{v} \cdot \mathbf{b} = 0$. We can think of $\mathbf{a}$ as a unit tangent vector to $M_1$ at $\mathbf{u}$, and $\mathbf{b}$ as a unit tangent vector to $M_2$ at $r \mathbf{v}$.  If we map $M_1$ and $M_2$ into a common copy of $\mathbb{R}^3$  so that the images of $M_1$ and $M_2$ are tangent to each other at the images of $\mathbf{u}$ and $r \mathbf{v}$, then the images of $\mathbf{a}$ and $\mathbf{b}$ coincide with each other on the common tangent plane.  Of course, we can choose another pair of unit tangent vectors to represent the same state: $\mathbf{u} \times \mathbf{a}$ and $-\mathbf{v} \times \mathbf{b}$ work nicely (if we assume that the two sphere images lie in each other's exteriors, as opposed to one sphere image nested inside the other).  For each point of contact $\left( \mathbf{u}, r \mathbf{v} \right)$, we can define an equivalence relation: two pairs of unit tangent vectors $\left( \mathbf{a}_1, \mathbf{b}_1 \right)$ and $\left( \mathbf{a}_2, \mathbf{b}_2 \right)$ define the same state at $\left( \mathbf{u}, r \mathbf{v} \right)$ if and only if both $\left<\mathbf{a}_1, \mathbf{a}_2 \right> = \left<\mathbf{b}_1, \mathbf{b}_2 \right>$ and $\left<\mathbf{a}_1, \mathbf{u} \times \mathbf{a}_2 \right> = -\left<\mathbf{b}_1, \mathbf{v} \times \mathbf{b}_2 \right>$.  We can define $\bar{M}$ to be the six dimensional manifold consisting on all $4$-tuples $\left( \mathbf{u}, \mathbf{v}, \mathbf{a}, \mathbf{b} \right)$, and we can define $M$, our rolling ball system, to be the five dimensional manifold consisting of all equivalence classes, or states, of those $4$-tuples.

The manifold $M$ inherits a natural sub-Riemmanian structure with growth vector $(2,3,5)$ from the Riemannian geometry of the spheres.  We interpret the ``no slipping, no twisting'' restriction on the motion of the two sphere to mean that a horizontal path $\gamma$ in $M$ is completely determined by motion along either sphere.  Specifically, for a $4$-tuple $\left( \mathbf{u}, \mathbf{v}, \mathbf{a}, \mathbf{b} \right)$, a unit speed path along the great circle on $M_1$ in the direction of $\mathbf{a} \cos \theta + (\mathbf{u} \times \mathbf{a}) \sin \theta$ will force unit speed motion along the great circle on $M_2$ in the direction of $\mathbf{b} \cos \theta - (\mathbf{v} \times \mathbf{b}) \sin \theta$.  This motion will include the parallel transport of $\mathbf{a}$ and $\mathbf{b}$ induced by the Levi-Civita connections on $M_1$ and $M_2$.  Altogether, a horizontal direction on $M$ will be a linear combination of the vectors $X_1 = \left( \mathbf{a}, \mathbf{b}, -\mathbf{u}, \displaystyle -\frac{\mathbf{v}}{r} \right)$ and $X_2 = \left( \mathbf{u} \times \mathbf{a}, -\mathbf{v} \times \mathbf{b}, \mathbf{0}, \mathbf{0} \right)$, which together form an orthonormal frame.  The distribution formed by these two vector fields is compatible with the equivalence relation defined above and therefore defines a two-dimensional horizontal distribution on $M$.

One of the nice features of $M$ is that it inherits a great deal of symmetry from $M_1$ and $M_2$.  Specifically, if $g_1, g_2 \in SO_3$, then the map $\left( \mathbf{u}, \mathbf{v}, \mathbf{a}, \mathbf{b} \right) \mapsto \left( g_1 \mathbf{u}, g_2 \mathbf{v}, g_1 \mathbf{a}, g_2 \mathbf{b} \right)$ is an isomorphism on $\bar{M}$, and therefore on $M$.  With this observation in mind, it suffices to describe the horizontal geodesics from any state $\mathbf{p} = \left[ \left( \mathbf{u}_0, \mathbf{v}_0, \mathbf{a}_0, \mathbf{b}_0 \right) \right]$ in $M$ to the state \begin{equation*} \mathbf{n} = \left[ \left( (0,0,1), (0,0,1), (1,0,0), (1,0,0) \right) \right] \end{equation*}the state where the spheres are in contact at their respective north poles, and the unit tangent vectors at these north poles both point towards the positive $x$-axis.  We can define a coordinate system $(w_1, w_2, w_3, w_4, w_5)$ for $M$ centered at $\mathbf{n}$ as follows: if $\left( \mathbf{u}, \mathbf{v}, \mathbf{a}, \mathbf{b} \right)$ is a representative of a state near $q$, then $w_1 = u_1$, $w_2 = u_2$, $w_3 = v_1$, $w_4 = v_2$, and \begin{equation*} w_5 = \left(a_1 u_3 - a_3 u_1\right) b_2+ \left( b_1 u_3 - b_3 u_2 \right) a_2 \end{equation*}  It is easy to show that the above coordinates do not depend on the representative used for the state.

Let $\gamma(t)$ be an absolutely continuous horizontal path on $M$ with a unit speed parametrization such that $\gamma(0) = \mathbf{p}$ and there exists $T$ such that $\gamma(T) = \mathbf{n}$.  Assume that $\gamma(t)$ extends past $\mathbf{n}$ to some point $\mathbf{q}$, and that $\gamma(t)$ is a length minimizing geodesic from $\mathbf{p}$ to $\mathbf{q}$.

Arbitrarily fix five times $t_1, t_2, t_3, t_4, t_5$ on $(0, T)$.  We are going to perturb $\gamma$ by bending its image at each of these five times.  Specifically, let $R_\mathbf{u}^\theta \in SO_3$ be the counterclockwise rotation about the axis passing through the unit vector $\mathbf{u}$ by an angle $\theta$ (here, the rotation is counterclockwise as viewed with $\mathbf{u}$ pointing towards us).  Let $\mathbf{\alpha} = (\alpha_1, \alpha_2, \alpha_3, \alpha_4, \alpha_5)$ be a 5-tuple of angles.  Define the the path $\gamma^\mathbf{\alpha}(t) = \left[ \left( \mathbf{u}^\alpha(t), \mathbf{v}^\alpha(t), \mathbf{a}^\alpha(t), \mathbf{b}^\alpha(t) \right) \right]$ as follows \begin{equation*} \mathbf{u}^\mathbf{\alpha}(t) = \begin{cases} \mathbf{u}(t) & 0 \leq t \leq t_1 \\ R_{\mathbf{u}(t_1)}^{\alpha_1} \mathbf{u}(t) & t_1 < t \leq t_2 \\ R_{\mathbf{u}(t_1)}^{\alpha_1} R_{\mathbf{u}(t_2)}^{\alpha_2} \mathbf{u}(t) & t_2 < t \leq t_3 \\ R_{\mathbf{u}(t_1)}^{\alpha_1} R_{\mathbf{u}(t_2)}^{\alpha_2} R_{\mathbf{u}(t_3)}^{\alpha_3} \mathbf{u}(t) & t_3 < t \leq t_4 \\ R_{\mathbf{u}(t_1)}^{\alpha_1} R_{\mathbf{u}(t_2)}^{\alpha_2} R_{\mathbf{u}(t_3)}^{\alpha_3} R_{\mathbf{u}(t_4)}^{\alpha_4} \mathbf{u}(t) & t_4 < t \leq t_5 \\ R_{\mathbf{u}(t_1)}^{\alpha_1} R_{\mathbf{u}(t_2)}^{\alpha_2} R_{\mathbf{u}(t_3)}^{\alpha_3} R_{\mathbf{u}(t_4)}^{\alpha_4} R_{\mathbf{u}(t_5)}^{\alpha_5} \mathbf{u}(t) & t > t_5   \end{cases} \end{equation*} \begin{equation*} \mathbf{v}^\mathbf{\alpha}(t) = \begin{cases} \mathbf{v}(t) & 0 \leq t \leq t_1 \\ R_{\mathbf{v}(t_1)}^{-\alpha_1} \mathbf{v}(t) & t_1 < t \leq t_2 \\ R_{\mathbf{v}(t_1)}^{-\alpha_1} R_{\mathbf{v}(t_2)}^{-\alpha_2} \mathbf{v}(t) & t_2 < t \leq t_3 \\ R_{\mathbf{v}(t_1)}^{-\alpha_1} R_{\mathbf{v}(t_2)}^{-\alpha_2} R_{\mathbf{v}(t_3)}^{-\alpha_3} \mathbf{v}(t) & t_3 < t \leq t_4 \\ R_{\mathbf{v}(t_1)}^{-\alpha_1} R_{\mathbf{v}(t_2)}^{-\alpha_2} R_{\mathbf{v}(t_3)}^{-\alpha_3} R_{\mathbf{v}(t_4)}^{-\alpha_4} \mathbf{v}(t) & t_4 < t \leq t_5 \\ R_{\mathbf{v}(t_1)}^{-\alpha_1} R_{\mathbf{v}(t_2)}^{-\alpha_2} R_{\mathbf{v}(t_3)}^{-\alpha_3} R_{\mathbf{v}(t_4)}^{-\alpha_4} R_{\mathbf{v}(t_5)}^{-\alpha_5} \mathbf{v}(t) & t > t_5   \end{cases} \end{equation*} \begin{equation*} \mathbf{a}^\mathbf{\alpha}(t) = \begin{cases} \mathbf{a}(t) & 0 \leq t \leq t_1 \\ R_{\mathbf{u}(t_1)}^{\alpha_1} \mathbf{a}(t) & t_1 < t \leq t_2 \\ R_{\mathbf{u}(t_1)}^{\alpha_1} R_{\mathbf{u}(t_2)}^{\alpha_2} \mathbf{a}(t) & t_2 < t \leq t_3 \\ R_{\mathbf{u}(t_1)}^{\alpha_1} R_{\mathbf{u}(t_2)}^{\alpha_2} R_{\mathbf{u}(t_3)}^{\alpha_3} \mathbf{a}(t) & t_3 < t \leq t_4 \\ R_{\mathbf{u}(t_1)}^{\alpha_1} R_{\mathbf{u}(t_2)}^{\alpha_2} R_{\mathbf{u}(t_3)}^{\alpha_3} R_{\mathbf{u}(t_4)}^{\alpha_4} \mathbf{a}(t) & t_4 < t \leq t_5 \\ R_{\mathbf{u}(t_1)}^{\alpha_1} R_{\mathbf{u}(t_2)}^{\alpha_2} R_{\mathbf{u}(t_3)}^{\alpha_3} R_{\mathbf{u}(t_4)}^{\alpha_4} R_{\mathbf{u}(t_5)}^{\alpha_5} \mathbf{a}(t) & t > t_5   \end{cases} \end{equation*} \begin{equation*} \mathbf{b}^\mathbf{\alpha}(t) = \begin{cases} \mathbf{b}(t) & 0 \leq t \leq t_1 \\ R_{\mathbf{v}(t_1)}^{-\alpha_1} \mathbf{b}(t) & t_1 < t \leq t_2 \\ R_{\mathbf{v}(t_1)}^{-\alpha_1} R_{\mathbf{v}(t_2)}^{-\alpha_2} \mathbf{b}(t) & t_2 < t \leq t_3 \\ R_{\mathbf{v}(t_1)}^{-\alpha_1} R_{\mathbf{v}(t_2)}^{-\alpha_2} R_{\mathbf{v}(t_3)}^{-\alpha_3} \mathbf{b}(t) & t_3 < t \leq t_4 \\ R_{\mathbf{v}(t_1)}^{-\alpha_1} R_{\mathbf{v}(t_2)}^{-\alpha_2} R_{\mathbf{v}(t_3)}^{-\alpha_3} R_{\mathbf{v}(t_4)}^{-\alpha_4} \mathbf{b}(t) & t_4 < t \leq t_5 \\ R_{\mathbf{v}(t_1)}^{-\alpha_1} R_{\mathbf{v}(t_2)}^{-\alpha_2} R_{\mathbf{v}(t_3)}^{-\alpha_3} R_{\mathbf{v}(t_4)}^{-\alpha_4} R_{\mathbf{v}(t_5)}^{-\alpha_5} \mathbf{b}(t) & t > t_5   \end{cases} \end{equation*} Some explanation is in order.  Consider $t > t_5$.  For such a time, we first create a counterclockwise bend in the image of $\gamma(t)$ on $M_1$ of angle $\alpha_5$ about the axis passing through $\mathbf{u}(t_5)$, and, simultaneously, a counterclockwise bend in the image on $M_2$ of angle $-\alpha_5$ about the axis passing through $\mathbf{v}(t_5)$.  At the same time, we also apply these rotations to the unit tangent vectors $\mathbf{a}(t)$ and $\mathbf{b}(t)$, respectively.  These rotations create discontinuities at $t_5$ for the functions $\mathbf{a}^\alpha(t)$ and $\mathbf{b}^\alpha(t)$, but the path through $M$, the space of equivalence classes, remains continuous.  We then repeat the same bending process at $t_4$, then at $t_3$, then at $t_2$, and finally at $t_1$, creating a new absolutely continuous path $\gamma^\mathbf{\alpha}(t)$ through $M$ whose images on $M_1$ and $M_2$ are bent at five points.

If all $\alpha_i$ is sufficiently near $0$, then $\gamma^\mathbf{\alpha}(T)$ will be a state in $M$ near $\mathbf{n}$ to use the coordinate system we defined above.  Let $\omega_j(\mathbf{\alpha}) = w_j \left(\gamma^\mathbf{\alpha}(T) \right)$, the $j$th coordinate of the image of $\gamma^\mathbf{\alpha}$ at time $T$ under the coordinate system defined above.  The map \begin{equation*} \omega \left( \alpha_1, \alpha_2, \alpha_3, \alpha_4, \alpha_5 \right) = \left( \omega_1(\mathbf{\alpha}), \omega_2(\mathbf{\alpha}), \omega_3(\mathbf{\alpha}), \omega_4(\mathbf{\alpha}), \omega_5(\mathbf{\alpha}) \right) \end{equation*} is clearly smooth, and if its derivative \begin{equation*} \begin{bmatrix} \frac {\partial \omega_1} {\partial \alpha_1} & \cdots & \frac {\partial \omega_5} {\partial \alpha_1} \\ \vdots & \ddots & \vdots \\ \frac {\partial \omega_1} {\partial \alpha_5} & \cdots & \frac {\partial \omega_5} {\partial \alpha_5} \end{bmatrix} \end{equation*} is non-singular, then $\omega$ maps open sets in $\mathbb{R}^5$ to open sets in $M$.  As a consequence, there exists $\bar{T} > T$ and $\alpha \in \mathbb{R}^5$ such that $\omega(\alpha)$ are the coordinates of $\gamma(\bar{T})$ under the coordinate system described above.  The arc length of $\gamma^\alpha$ from time $0$ to time $t$, however, is still $t$, since the construction of $\gamma^\alpha$ above does not change arc length.  Therefore $\gamma$ cannot be the length minimizing curve from $\mathbf{p}$ to $\gamma(\bar{T})$, since $\gamma^\alpha$ is shorter, but this leads to a contradiction.  Thus the derivative of $\omega$ must be singular at $\mathbf{\alpha} = (\alpha_1, \alpha_2, \alpha_3, \alpha_4, \alpha_5) = (0,0,0,0,0)$ if $\gamma$ is a length minimizing geodesic.

To take advantage of this observation, we examine the derivatives $\displaystyle \frac {\partial \omega_j} {\partial \alpha_i}$.  To compute $\displaystyle \frac {\partial \omega_j} {\partial \alpha_i}$, we fix all angles at $0$ except $\alpha_i$.  Thus computing these derivatives is reduced to understanding the actions of rotating counterclockwise about the axes passing through $\mathbf{u}^\mathbf(t_i)$ and $\mathbf{v}(t_i)$.  These actions will only depend on the locations of $\mathbf{u}(t_i)$ and $\mathbf{v}(t_i)$.  Specifically, we compute \begin{equation*} \frac {\dif \mathbf{u}^\mathbf{\alpha}(t)} {\dif \alpha_i} = (u_2(t_i), -u_1(t_i), 0) \qquad \frac {\dif \mathbf{v}^\mathbf{\alpha}(t)} {\dif \alpha_i} = (-rv_2(t_i), rv_1(t_i), 0) \end{equation*} \begin{equation*} \frac {\dif \mathbf{a}^\mathbf{\alpha}(t)} {\dif \alpha_i} = (0, u_3(t_i), u_2(t_i)) \qquad \frac {\dif \mathbf{b}^\mathbf{\alpha}(t)} {\dif \alpha_i} = (0, -v_3(t_i), -v_2(t_i)) \end{equation*}  Immediately this gives us derivatives for the first four coordinates: \begin{equation*} \frac {\partial \omega_1} {\partial \alpha_i} = u_2(t_i) \qquad \frac {\partial \omega_2} {\partial \alpha_i} = -u_1(t_i) \end{equation*} \begin{equation*} \frac {\partial \omega_3} {\partial \alpha_i} = -rv_2(t_i) \qquad \frac {\partial \omega_4} {\partial \alpha_i} = rv_1(t_i) \end{equation*} The derivative of the last coordinate $\omega_5$ is only slightly more difficult: it requires just an application of the product rule and the fact the $\gamma^\mathbf{\alpha}(t)$ passes through $\mathbf{n}$ at time $T$: \begin{equation*} \frac {\partial \omega_5} {\partial \alpha_i} = u_3(t_i) - v_3(t_i) \end{equation*} Using the calculations above, we can fill in the matrix\begin{equation*} \begin{bmatrix} \frac {\partial \omega_1} {\partial \alpha_1} & \cdots & \frac {\partial \omega_5} {\partial \alpha_1} \\ \vdots & \ddots & \vdots \\ \frac {\partial \omega_1} {\partial \alpha_5} & \cdots & \frac {\partial \omega_5} {\partial \alpha_5} \end{bmatrix} = \begin{bmatrix} u_2(t_1) & -u_1(t_1) & -r v_2(t_1) & r v_1(t_1) & u_3(t_1) - v_3(t_1) \\ u_2(t_2) & -u_1(t_2) & -r v_2(t_2) & r v_1(t_2) & u_3(t_2) - v_3(t_2) \\ u_2(t_3) & -u_1(t_3) & -r v_2(t_3) & r v_1(t_3) & u_3(t_3) - v_3(t_3) \\ u_2(t_4) & -u_1(t_4) & -r v_2(t_4) & r v_1(t_4) & u_3(t_4) - v_3(t_4)  \\ u_2(t_5) & -u_1(t_5) & -r v_2(t_5) & r v_1(t_5) & u_3(t_5) - v_3(t_5)  \end{bmatrix} \end{equation*} As we argued before, the above matrix must be singular.  Therefore, there must exist constants $k_1$, $k_2$, $k_3$, $k_4$, and $k_5$, not all equal to $0$, such that, for all $i = 1, 2, 3, 4, 5$, we have that \begin{equation*} k_1 u_1(t_i) + k_2 u_2(t_i) + k_3 v_1(t_i) + k_4 v_2(t_i) + k_5 \left(u_3(t_i) - v_3(t_i) \right) = 0 \end{equation*}
The constants $k_i$ are determined by four of the times $t_i$, and replacing the remaining time with any other $t \in [0, T]$ yields the same equation.  Thus, for any time $t \in [0, T]$, we must have \begin{equation*} k_1 u_1(t) + k_2 u_2(t) + k_3 v_1(t) + k_4 v_2(t) + k_5 \left(u_3(t) - v_3(t) \right) = 0 \end{equation*}  Because of the homogeneity of $M$, there is nothing unique about $\mathbf{n}$, so a similar plane equation holds for length minimizing geodesics passing through any other point in $M$.  Thus, if $\gamma$ is a length minimizing geodesic, there exists a 5-plane $P$ in $\mathbb{R}^3 \times \mathbb{R}^3$ such that all points of contact $\left( \mathbf{u}(t), \mathbf{v}(t) \right) \in \mathbb{R}^3 \times \mathbb{R}^3$ for $\gamma$ lie on $P$.

\end{document}